\documentclass{SCIS2018}
\begin{document}

\ArticleType{LETTER}
\Year{2018}
\Month{}
\Vol{61}
\No{}
\DOI{}
\ArtNo{}
\ReceiveDate{}
\ReviseDate{}
\AcceptDate{}
\OnlineDate{}


\title{Synthesis of model predictive control based on data-driven learning}
{Synthesis of model predictive control based on data-driven learning}

\author[1, 2, 3]{Yuanqiang ZHOU}{}
\author[2, 3]{Dewei LI}{{dwli@sjtu.edu.cn}}
\author[2, 3]{Yugeng XI}{}
\author[1]{Zhongxue GAN}{}

\AuthorMark{Zhou Y}

\AuthorCitation{Zhou Y, Li D, Xi Y, Gan Z}


\address[1]{ENN Science \& Technology Development Co.,Ltd, Langfang, Hebei, {\rm 065001}, China}
\address[2]{Department of Automation, Shanghai Jiao Tong University, Shanghai {\rm 200240}, China}
\address[3]{Key Laboratory of System Control and Information Processing, Ministry of Education of China, Shanghai {\rm 200240}, China}

\maketitle

\begin{multicols}{2}
\deareditor
Model predictive control (MPC) is a practically effective and attractive approach in the field of industrial processes~\cite{1} owing to its excellent ability to handle constraints, nonlinearity, and performance/cost trade-offs. The core of all model-based predictive algorithms is to use ``open-loop optimal control" instead of ``closed-loop optimal control" within a moving horizon~\cite{2}. It is assumed in this letter that the reader is familiar with MPC as a control design methodology. 

Because the dynamic model of a system predicting its evolution is usually inaccurate, the actual behaviors may deviate significantly from the predicted ones. Thus, acquiring accurate knowledge of the physical model is essential to ensure satisfactory performance of MPC controllers. Owing to the well-developed information technology, copious amounts of measurable process data can be easily collected, and such data can then be employed to predict and assess system behaviors and make control decisions, especially for the establishment and development of learning MPC.

For the application of MPC design in on-line regulation or tracking control problems, several studies have attempted to develop an accurate model, and realize adequate uncertainty description of linear or non-linear plants of the processes~\cite{3,4,5}.
In this work, we employ the data-driven learning technique specified in~\cite{6} to iteratively approximate the dynamical parameters, without requiring \emph{a priori} knowledge of system matrices. The proposed MPC approach can predict and optimize the future behaviors using multi-order derivatives of control input as decision variables. 
Because the proposed algorithm can obtain a linear system model at each sampling, it can adapt to the actual dynamics of time-varying or nonlinear plants. This methodology can serve as a data-driven identification tool to study adaptive optimal control problems for unknown complex systems. 

\lettersection{Problem Formulation}
In this work, we consider a continuous-time industrial process given by
\begin{align}\label{eq:1}
{\dot x}(t)  = &  A x(t) + B u(t) \notag \\
 \buildrel \vartriangle \over = &  \mathcal{H} \left(x(t), u(t) \right) {\Theta} 
\end{align} 
where $t \ge {t_0}$, ${x} \in \mathbb{R}^{n}$, and $ u\in \mathbb{R}^{m}$ are the system states and input, respectively. $\mathcal{H} (\cdot, \cdot): \mathbb{R}^n \times \mathbb{R}^m \rightarrow \mathbb{R}^{n \times ({n^2 + mn})} $ is defined as $\mathcal{H}(x, u) \buildrel \vartriangle \over =  \left[(x \otimes I_n)^{\textrm{T}} ~  ( u \otimes I_n  )^{\textrm{T}} \right]$, where $\otimes$ denotes the Kronecker product. $\Theta$ denotes the vector of the system parameters given by $\Theta \buildrel \vartriangle \over =\left[ \textrm{vec}(A)^{\textrm{T}}  ~ \textrm{vec}(B)^{\textrm{T}} \right]^{\textrm{T}}   \in \mathbb{R}^{n^2 + nm}$, where $A \in \mathbb{R}^{n \times n}$ is the system matrix, $B \in \mathbb{R}^{n \times m}$ is the input matrix, and $\mathrm{vec}(\cdot)$ denotes the vectorization operator, that is, $\textrm{vec}(P) = \left[p_1^\textrm{T}, \ldots, p_m^\textrm{T}\right]^\textrm{T}$, where $p_i \in \mathbb{R}^n$ is the $i$th column of a matrix $P \in \mathbb{R}^{n \times m}$. We assume that $(A,B)$ is controllable and $(A, C)$ is observable.

In this work, we consider the following input constraint: $u \in \mathcal{U} \subset \mathbb{R}^m$, where $\mathcal{U}$ denotes a nonempty compact convex set and contains the origin as its interior point. In this case, since $\Theta$ is unknown, the primary objective of this work is to design a data-driven MPC formulation to obtain an open-loop optimal control policy that tracks a given reference $x_{\textrm{d}}$ and, at each sampling time $t_k, k=1,2, \ldots$, minimizes the following cost function
\begin{align}\label{eq:2}
J(x(t_k), \hat u_k(\cdot)) ) \!=\! & \int_{t_k }^{ t_k+ T} \big( \|e (\tau)\|_Q^2  +\| \hat u_k (\tau) ) \|_R^2 \big) \textrm{d} \tau \notag \\
& + \varPhi (e (t_k + T)), 
\end{align}
where $e(\cdot) = x(\cdot) -x_{\textrm{d}}(\cdot)$ denotes the error, $\varPhi (\cdot)$ denotes the terminal cost, and $Q  = Q^T \! \succ \! 0$ and $R = R^T \! \succeq \! 0$ are the symmetric weighting matrices.

\lettersection{Methodology}
To facilitate MPC design, at time $t =t_k$, the states $x$ and the parameter $\Theta$ of the predicted model over the moving horizon $[t, t + T]$ are both learned from the input--output measurements, using a data-driven learning technique. 
In this work, we consider two situations.
\begin{itemize}
\item All the states and input information are available to us. Then, by rearranging (1), we have the linear error system in the form
	\begin{align}\label{eq:3}
	\mathcal{F} (t)  = \mathcal{G} (t) \hat \Theta, \quad  \forall t  \in \mathbb{R}_{\ge 0}, 
	\end{align}
	where $\hat \Theta$ is an estimate of the unknown parameter $\Theta$; the matrices $\mathcal{F} (\cdot): \mathbb{R}_{\ge 0} \rightarrow \mathbb{R}^n $ and $\mathcal{G} (\cdot): \mathbb{R}_{\ge 0} \rightarrow \mathbb{R}^{n \times (n^2 + mn)} $ are defined as
	\begin{align*} 
	\mathcal{F} (t) = & \begin{cases}   x(t) - x(t - \delta), \quad t \in [\delta, \infty),  \\
	0, \qquad \qquad \qquad \quad t<\delta, 
	\end{cases}  \\
	\mathcal{G} (t) = & \left[(\Xi_x(t)  \otimes I_n)^\textrm{T} ~~ (\Xi_u (t)  \otimes I_n)^\textrm{T} \right], 
	\end{align*}
	where $\delta$ denotes the sampling period; the vectors $\Xi_x (\cdot) : \mathbb{R}_{\ge 0} \rightarrow \mathbb{R}^n$ and $\Xi_u (\cdot): \mathbb{R}_{\ge 0} \rightarrow \mathbb{R}^m $ are defined as
	\begin{align*}
    \Xi_x(t) = & \begin{cases}   \int_{t - \delta}^{t } x(\tau ) \textrm{d} \tau, \quad t \in [\delta, \infty),  \\
	0, \qquad \qquad \quad \quad t<\delta,
	\end{cases}    \\
	\Xi_u (t) = & \begin{cases}   \int_{t - \delta}^{t } u(\tau )  \textrm{d} \tau, \quad  t \in [\delta, \infty),  \\
	0, \qquad \qquad \quad \quad  t<\delta.
	\end{cases}  
	\end{align*}
	
\item Only partial states and input information are available; we assume the available states as the first $q = n/2 < n$ components of the states and denote them as $\xi \in \mathbb{R}^q$. We assume that the pair $(A, B)$ has the form
\begin{align}\label{eq:4}
A= \begin{bmatrix}
0_{q \times q} &  I_{q} \\
A_1 & A_2
\end{bmatrix}, \quad B= \begin{bmatrix}
0_{q \times m}  \\
B_1
\end{bmatrix}.
\end{align}
Then, the linear error system is given by
\begin{align}\label{eq:5}
\mathcal{F}_1 (t)  = \mathcal{G}_1 (t) \hat \Theta_1, \quad  \forall t  \in \mathbb{R}_{\ge 0}, 
\end{align}
where $\hat \Theta_1$ is an estimate of the unknown parameter $\Theta_1 =\left[ \textrm{vec}(A_1)^{\textrm{T}}  ~\textrm{vec}(A_2)^{\textrm{T}}  ~ \textrm{vec}(B_1)^{\textrm{T}} \right]^{\textrm{T}} \in \mathbb{R}^{2 q^2 + mq}$; $\mathcal{F}_1 (\cdot): \mathbb{R}_{\ge 0} \rightarrow \mathbb{R}^q $ and $\mathcal{G}_1 (\cdot): \mathbb{R}_{\ge 0} \rightarrow \mathbb{R}^{q \times (2 q^2 + mq)} $ are defined as
\begin{align*} 
\mathcal{F}_1 (t) = & \begin{cases}   \xi(t \!-\delta_2 \!-\! \delta_1) \!-\! \xi(t - \delta_1) \!+\! \xi(t)  \\
\qquad   -\xi(t-\delta_2), \quad t \in [\delta_1 + \delta_2, \infty),  \\
0, \quad \qquad \qquad \qquad t<\delta_1 + \delta_2,
\end{cases}  \\
\mathcal{G}_1 (t) = & \left[(\Xi_p(t)  \otimes I_n)^\textrm{T} ~~(\Xi_v(t)  \otimes I_n)^\textrm{T} \right. \notag \\ & \qquad \qquad  \qquad  \qquad    \left. (\Xi_u^1(t)  \otimes I_n)^\textrm{T} \right],
\end{align*}
where $\delta_1$ and $ \delta_2$ ($\delta_1 \neq \delta_2$) denote the different periods; $ \Xi_p(\cdot) : \mathbb{R}_{\ge 0} \rightarrow \mathbb{R}^q$, $\Xi_v (\cdot): \mathbb{R}_{\ge 0} \rightarrow \mathbb{R}^q$, and $\Xi_u^1 (\cdot): \mathbb{R}_{\ge 0} \rightarrow \mathbb{R}^m$ are defined as
\begin{align*}
\Xi_p(t) = & \begin{cases}  \int_{t - \delta_2}^{t }  \int_{\tau  - \delta_1}^{ \tau } \xi (\varsigma ) \textrm{d}\varsigma \textrm{d} \tau, ~ t \in [\delta_1 + \delta_2, \infty),  \\
0, \qquad \qquad \quad \quad t < \delta_1 + \delta_2, 
\end{cases}  \\
\Xi_v(t) = & \begin{cases}   \int_{t \!-\! \delta_2}^{t } \xi (\tau ) \textrm{d} \tau -  \int_{t \!-\! \delta_1 \!-\! \delta_2 }^{t \!-\!\delta_1} \xi(\tau )  \textrm{d} \tau,  \\
\qquad \qquad  \qquad \qquad   t \in [\delta_1 + \delta_2, \infty),  \\
0, \qquad \qquad \quad \qquad  t<\delta_1 + \delta_2,
\end{cases}  \\
\Xi_u^1 (t) = & \begin{cases}    \int_{t - \delta_2}^{t }  \int_{\tau  - \delta_1}^{ \tau } u (\varsigma ) \textrm{d} \varsigma \textrm{d} \tau, ~ t \in [\delta_1 + \delta_2, \infty),  \\
0, \qquad \qquad \quad \quad \qquad \quad  t < \delta_1 + \delta_2.
\end{cases} 
\end{align*}
\end{itemize}

Further, from (3) and using the measurements, for a positive integer $l \le k$, we define the vector $\varGamma_k \in \mathbb{R}^{l n} $ and matrix $\Psi_k \in \mathbb{R}^{l n \times (n^2 + mn)}$ such that
\begin{align*}
\varGamma_k  \buildrel \vartriangle \over =  & \Big[ \mathcal{F}^\mathrm{T} (t_0), \mathcal{F}^\mathrm{T} (t_1), \ldots, \mathcal{F}^\mathrm{T} (t_l)  \Big]^\mathrm{T} ,   \\
\Psi_k  \buildrel \vartriangle \over = & \Big[ \mathcal{G}^\mathrm{T} (t_0), \mathcal{G}^\mathrm{T} (t_1), \ldots, \mathcal{G}^\mathrm{T} (t_l)  \Big]^\mathrm{T} , 
\end{align*}
where $0 \le t_0 < t_1 < \cdots < t_l$ and $t_i = i \delta, i = 0, 1, \ldots, l$. 
Then, (3) implies the linear equation
\begin{align} \label{eq:6}
\varGamma_k  =  \Psi_k \hat \Theta, 
\end{align}
Notice that if $\Psi_k$ has full column rank, (6) can be directly solved as
\begin{align} \label{eq:7}
\hat \Theta  = \left( \Psi_k^\mathrm{T} \Psi_k \right)^{-1} \Psi_k^\mathrm{T} \varGamma_k. 
\end{align}
Similarly, for (5), we let $\delta_1 = \delta$ and $\delta_2 = 2 \delta$; thus, we obtain the same results for $\hat \Theta_1$. 
To guarantee $\text{rank}( \Psi_k^\mathrm{T} \Psi_k ) = n^2 + nm$, we let the states and inputs collected over a sufficiently large number of data samples be $l \gg n^2 + nm$. In practice, we assume that there exists a nominal control input $u = -K_0 x$, where $K_0$ denotes a stabilizing feedback gain matrix, such that $\varGamma_k$ and $ \Psi_k$ in (6) can be implemented using $2l$ integrators to collect information about the states and inputs.  Using $\mathfrak{D}_k \buildrel \vartriangle \over =   \bigcup\limits_{i = 0}^l { \{\mathcal{F} (t_i), \mathcal{G} (t_i) \} }$ by (7), we have $\hat \Theta  = \left[ \mathrm{vec} (\hat A)^\mathrm{T}  ~\mathrm{vec} (\hat B)^\mathrm{T} \right]^\mathrm{T}$, which can be used to predict and optimize the future behaviors over a finite horizon $[0, T]$. 

First, let us consider the MPC formulation. As mentioned in~\cite{4}, the majority (if not all) of existing formulations consider only $u(\cdot)$ as the decision variable. For (1), we can extend it to a higher order derivative of $u(\cdot)$, that is,
\begin{align}
{\mathbf{u}}_k (t)  \buildrel \vartriangle \over = \Big[ \hat u_k^\textrm{T}(t), (\hat u_k^{[1]})^\textrm{T}(t), \ldots, (\hat u_k^{[r]})^\textrm{T}(t) \Big]^\textrm{T}, 
\end{align}
with some control order $r \in \mathbb{N}_{+}$ larger than $\rho \ge 1$, where $\rho$ denotes the input relative degree of (1). This will improve the efficacy of our learning MPC, and the first term of $\mathbf{u}_k (\cdot)$ in (8) is the control input $\hat u_k (\cdot)$ that is to be optimized in (2). Then, we let $ \tilde B = \textrm{vec}^{-1}(\hat B)$,  $ \tilde A = \textrm{vec}^{-1}(\hat A)$, where $\textrm{vec}^{-1}(\cdot)$ denotes the inverse operation of $\textrm{vec}(\cdot)$, and define the following matrices: 
\begin{align*}
	\mathcal{A}_1  \buildrel \vartriangle \over = &  \begin{bmatrix}  I & \tilde A^\textrm{T}   & \cdots &   (\tilde A^{\rho - 1})^\textrm{T}   \end{bmatrix}^\textrm{T} ,  \\    
	\mathcal{A}_2  \buildrel \vartriangle \over = &  \begin{bmatrix}   (\tilde A^{\rho })^\textrm{T}  & (\tilde A^{\rho +1 })^\textrm{T}  & \cdots &   (\tilde A^{ r })^\textrm{T}   \end{bmatrix}^\textrm{T}, \mathrm{and}\\
	\mathcal{B}  \buildrel \vartriangle \over = &  \begin{bmatrix}   \tilde A^{\rho-1 } \tilde B & 0 &  \cdots & 0 \\ \tilde A^{\rho }  \tilde B&   \tilde A^{\rho-1 } \tilde B  & \cdots & 0  \\ \vdots & \vdots & \ddots & \vdots \\   {\tilde A}^r  \tilde B &  {\tilde  A}^{r -1}  \tilde B & \cdots & { \tilde A}^{ \rho -1 } \tilde B  \end{bmatrix}.
\end{align*}
At time instant $t = t_k$, the MPC formulation can be given by Eq. (9) (see Appendix A for details). 
\begin{align}
& \qquad  {\mathbf{u}}_k^\star(\cdot) =\mathop{ \arg \min }\limits_{{{ {\mathbf{u}}_k }}(\cdot)} {J}({{x}}({t_k}), {{\mathbf{u}}_k } (\cdot)) \label{eq:25}   \\
&\textrm{s.t.} \quad x(t_k+ \tau )  = \begin{bmatrix}  T_1 (\tau) & T_2 (\tau) \end{bmatrix}   \begin{bmatrix}
X_1 \\  X_2 \end{bmatrix}, \tau \in [0, T]  \notag  \\
& \qquad ~ X_1  = \mathcal{A}_1 x, \quad X_2 =  \mathcal{A}_2 x +  \mathcal{B} {\mathbf{u}}_k,  \notag \\
& \qquad ~ {\hat x}(t_k) = x(t_k), \quad {\hat u_k}(t) \in \mathcal{U},  \notag
\end{align}
where ${X}_1=\left[ x^\textrm{T}, (x ^{[1]})^\textrm{T} , \ldots,  (x ^{[\rho - 1 ]}) ^\textrm{T}  \right]^\textrm{T}$, ${X}_2=\left[  (x^{[\rho]} ) ^\textrm{T}, (x^ {[\rho + 1 ]})^\textrm{T} , \ldots,  (x^{[r]})^\textrm{T} \right]^\textrm{T}$, $ T_1 (\tau) =  \left[ 1, \tau, \ldots,  \frac{\tau^{\rho -1}}{(\rho - 1)!}  \right]$, $T_2(\tau) = \left[    \frac{\tau^{\rho}}{ \rho !}, \ldots,  \frac{\tau^{ r }}{ r!}  \right]$, $T_3(\tau) = \left[T_1(\tau), T_2(\tau) \right]$, $X_{1,\textrm{d}} = \left[ x_\textrm{d}^\textrm{T}, (x_\textrm{d}^{[1]})^\textrm{T} , \ldots ,  (x_\textrm{d}^{[\rho - 1 ]}) ^\textrm{T}  \right]$, 
$X_{2, \textrm{d}} =\left[  (x_\textrm{d}^{[\rho]} ) ^\textrm{T}, (x_\textrm{d}^ {[\rho + 1 ]})^\textrm{T} , \ldots ,  (x_\textrm{d}^{[r]})^\textrm{T} \right]$, and 
\begin{align*}
{J}  ({{x}}({t_k}), & {{\mathbf{u}}_k } (\cdot))  = \tilde X_1^\textrm{T} {\mathcal{T}}_{1,1} \tilde X_1 + 2 \tilde{X}_1^\textrm{T} {\mathcal{T}}_{1,2} \tilde {X}_2 \\
& +  \tilde{X}_2^\textrm{T} {\mathcal{T}}_{2,2} \tilde{X}_2 +   { {\mathbf{u}}_k} ^\textrm{T}  \mathcal{T}   { {\mathbf{u}}_k}   +  \varPhi (\tilde X_i (t_k + T)), 
\end{align*}
with $\tilde X_i = X_i -X_{i,\textrm{d}}$, $\Xi_i (\tau)= \sqrt{Q}  T_i(\tau)$, ${\mathcal{T}}_{i, j} = \int_{0}^{T} \Xi_i^\textrm{T} \Xi_j \textrm{d}\tau$, $i, j \in \{1, 2\}$, and ${\mathcal{T}} = \int_{0}^{T}  T_3^\textrm{T} R T_3 \textrm{d}\tau$.
This is a standard quadratic programming (QP) problem that can be solved by many available tools. In particular, for the case with box constraints,  with consideration of the possible model error with the data-driven method, we can handle the input constraints by the sub-optimal method in Appendix B.  We thus have the optimal control policy $\hat u_k^\star (t) = I_u   {\mathbf u}_k^\star$, where $I_u = [1, 0, \ldots, 0]_{1 \times (r+1)}$. We summarize this proposed approach as Algorithm 1 in  Appendix C. 

Furthermore, we consider the linear error system with the control policy $\hat u_k^\star (t)$ applied to (1). For the actual state trajectory, we have the continuous error as
\begin{align}\label{eq:22}
w(t) =   \mathcal{H} \left(x(t), \hat u_k^\star (t) \right) \Big( \hat{\Theta} -  \Theta  \Big),  
\end{align} 
where $ t \in [t_k, t_k + T]$. In Appendix D, we show that $w(t)$ is bounded and has the upper
bounded rate of change with time $t$; if $\tilde A = A$ and $\tilde B = B$, then $w(t) = 0, t\ge 0$; with the updated control policy $\hat u_k^\star(t)$ at each time $t = t_k$, $\lim_{t \rightarrow \infty}  w(t) =0 $, which implies the asymptotic stability of the closed-loop system. 

\lettersection{Simulation results}
An illustrative numerical example for two continuous stirred tank reactor (CSTR) systems is provided to validate the performance of the proposed approach. More details and discussions are presented in Appendix E. 

\Acknowledgements{This work was supported by the National Key Basic Research Special Foundation of China (2014CB249200) and the authors would like to thank Prof. Zhong-ping Jiang and his CAN Lab at Tandon School of Engineering, New York University, Brooklyn, NY, USA, for many inspirations and the help of this work. }

\Supplements{Appendix A, B, C, D and E.}


\end{multicols}
\end{document}